\documentclass[10pt]{amsart}
\usepackage{graphicx}

\newcommand{\N}{\mathbb{N}}
\newcommand{\R}{\mathbb{R}}

\newcommand{\be}{\begin{equation}}
\newcommand{\ee}{\end{equation}}

\newtheorem{theorem}{Theorem}[section]
\newtheorem{proposition}[theorem]{Proposition}
\newtheorem{lemma}[theorem]{Lemma}

\newtheorem{remark}[theorem]{Remark}
\newtheorem{definition}[theorem]{Definition}

\begin{document}

\thispagestyle{empty}

\title[Forward-backward parabolic problems]{Two--phase entropy solutions of forward--backward parabolic problems with unstable phase}

\author{Andrea Terracina}
\maketitle

\begin{abstract} 
This paper study the two--phase problem for the forward--backward parabolic equation with diffusion function of cubic type. Existence and uniqueness for these kind of problems were obtained in literature  in the case in which  the phases are both stable.  Here we consider the situation in which the unstable phase is taken in account, obtaining not trivial solution of the problem.   It is interesting to note that such solutions are given by solving generalized Abel's equations. \end{abstract}


\section{Introduction}\label{intro}

 In this paper we study  the following    {\it forward--backward}
parabolic problem: 
\be\label{B1}
 \left\{\begin{array}{lll} u_t= \phi(u)_{xx} & \hbox{ in } Q_T:=\R\times (0,T)\\
  u(x,0)=u_0(x) & \hbox{ in }\R\times\{0\},\end{array}
  \right.
 \ee
where   $\phi$ is a {\it nonmonotone}   function. This problem arises in different mathematical model: phase transition \cite{BS}, population dynamics \cite{Pa1}, \cite{Pa2}, oceanography \cite{BBDU}, image processing  \cite{PM}.
Obviously  problem \eqref{B1} is ill--posed
whenever $u$ takes values in the interval in which $\phi$ decreases.  
 
 In this paper we focus to the model of phase transition where the response function $\phi$ is of ``cubic type''.  More precisely we assume that $\phi\in Lip_{loc}(\R)$ and   
 $$\lim_{u\to\pm\infty}\phi(u)=\pm \infty.$$
  We suppose that $\phi$ has  a local minimum $A$ and a local maximum $B$ such that $A<B$. Let us denote with $c$, respectively $b$, the point in which the local minimum $A$, respectively maximum $B$, is achieved.  There are  three regions $(-\infty,b)$,   $(c,+\infty)$ and $(b,c)$;  the first two in which $\phi$ increase and the last one in which $\phi$ decrease.
 
 In the phase transition models, the function $u$ gives the phase fields, then  the  increasing intervals correspond to the stable phases and the interval  $(b,c)$ to the unstable, or metastable one.
 
 In this framework it is included also the piecewise linear case in which $\phi$  is given by 
 \begin{equation}\label{philin}
\phi(u) = \left\{\begin{aligned}
&\phi_1(u) \qquad & \textrm{for}\quad  & u \le b \\
&\phi_0(u)\qquad   & \textrm{for}\quad  & b < u < c\\
&\phi_2(u)\qquad  & \textrm{for}\quad   & u \ge d \, ,
\end{aligned}\right.
\end{equation}
where
\begin{equation*}
\phi_i(u):=\gamma_i\,u+\delta_i \, ,\,i=0,1,2; .
\end{equation*}
\begin{figure}[htbp]\label{funzione2}
\begin{center}
\includegraphics[width=9cm, height=5cm]{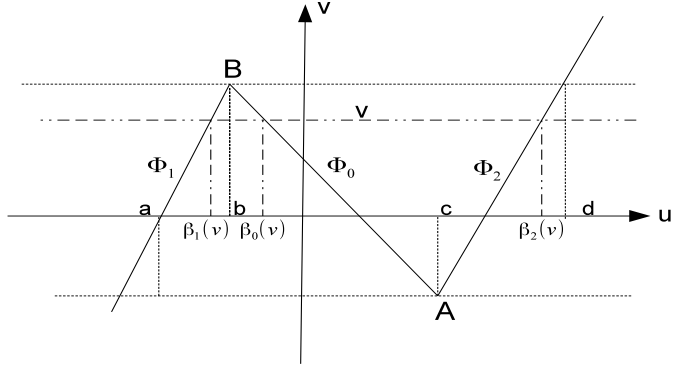}
\caption{The function  $\phi$. }
\end{center}
\end{figure}

Here $-\infty<b<c<\infty$,  $\alpha_i>0$, 
$\gamma_i\in\R$, $i=1,2$,
$A:=\phi_2(c)<\phi_1(b)=:B$.
In particular, 
$\phi_0(u):= \frac{A(u-b) -B(u-c)}{c-b}.$

 In \cite {H} (see also \cite{HN}), it was proved that  uniqueness does not hold for  problem (\ref{B1})  also in the class of  solutions that take value only in the two stable phases. 
 
Therefore, it is necessary to impose some stricter conditions in order to give a good formulation for the problem (\ref{B1}).  The idea is to introduce a proper regularization that comes from the physical phenomena of the original model. A classical approximation term is that introduced by the Cahn--Hilliard model that  describes the cost of the  inhomogeneities in phase transition.  On the other hands the mathematical description  of the physical phenomena is much more complicated to that  given by the Cahn--Hilliard equation and there   are some other terms to take in account (see e.g.  \cite {BFJ}, \cite {BS}, \cite{fife}, \cite{gurtin} \cite{BBDU}, \cite{BFG},   \cite{GG}, \cite{V}). As a matter of fact, it is possible to choose different type of regularizations  in which only some  phenomena are  highlighted.

In this case we refer to the following pseudoparabolic regularization 
\begin{equation} \label{E2}
 \left\{\begin{array}{lll} u_t= v_{xx}  & \hbox{ in } Q_T=\R\times (0,T)\\
 u(x,0)=u_0(x) & \hbox{ in } \R\times\{0\},\end{array}
  \right.\ee
where $v=\phi(u)+\epsilon u_t$, $\epsilon>0$.

\noindent The third order term  in the right hand   side of the differential equation in (\ref{E2})  is a viscosity term  related to   nonequilibium effects (see e.g. \cite{BFJ}, \cite{FJ}, \cite{gurtin}).

In \cite{NP} it was considered  the following Neumann boundary problem 

\begin{equation} \label{NP2}
 \left\{\begin{array}{lll} u_t= \Delta{(\phi(u)+\epsilon u_t)}  & \hbox{ in } Q_T=\Omega\times (0,T)\\
 \frac{\partial}{d\nu}(\phi(u)+\epsilon u_t)=0& \hbox{ in } Q_T=\partial\Omega\times (0,T)\\
  u(x,0)=u_0(x) & \hbox{ in } \Omega\times\{0\},\end{array}
  \right.\ee
where $\Omega\subset \R^n$ is a bounded ``regular'' domain. It was obtained  global existence and uniqueness for a large class of initial data. Moreover,  it was proved that the solutions satisfy some viscous entropy inequalities that are  crucial  to obtain a priori estimates that do not depend on the parameter $\epsilon$.

The singular limit  of problems \eqref{NP2}  was  analyzed  by Plotnikov (see \cite{Pl1}, \cite{Pl2}  \cite{Pl3}). Using this approach, it was obtained  an entropy formulation for the solution of the original forward--backward problem,  assuming that the physical solutions  of it  are that obtained when $\epsilon$ goes to $0^+$ as limit of solutions of problems (\ref{NP2}) (see e.g \cite{Pl2}, \cite{MTT}, \cite{T2}).  We do not enter in the details of the formulation given by Plotnikov (see Section \ref{2}). However it is important to mention that the original forward--backward equation is not more satisfied.   Just to give an idea, we  can extract subsequences $\{u_{\epsilon_n}\}$ of solutions of the viscous problems \eqref{NP2}  converging in the weak * topology of $L^{\infty}$ to a function $u$ and a  corresponding sequence $\{v_{\epsilon_n}\}$, defined by  $v_{\epsilon_n}=\phi(u_{\epsilon_n})+\epsilon_n (u_{\epsilon_n})_t$, such that, in general, converges to $v\neq\phi(u)$.
 Then, we have  $u_t=v_{xx}$  in a weak way. However, in general, the  original forward-backward equation does not hold. Plotnikov gave a characterization of the Young measure associated to the converging sequence $\{u_{\epsilon_n}\}$ proving that this is given by  a superposition of three Dirac measures. In particular there is the following relation between the functions $u$ and $v$:
\be \label {lambda}u=\lambda_1\beta_1(v)+\lambda_0\beta_0(v)+\lambda_2\beta_2(v),\ee
where $\lambda_i(x,t)\ge 0$,  $i=0,1,2$, $\displaystyle{\sum_{i=0}^2\lambda_i(x,t)=1}$ in $Q_T$, $\beta_i(s)$ correspond to the three monotone branches of the  graph  $s=\phi(u)$ (see Fig. \ref{funzione2}). 
From a physical point of view this fact could be interpreted as a superposition of different phases where $\lambda_i$ represents  the fraction of phase $i$.

\noindent Then, it could be proved that the equation $u_t=\phi(u)_{xx}$ is satisfied in the sense of the measure--valued solutions (see  \cite{Pl2}, \cite{MTT}). 

\noindent We can guess that  the complicate structure of the solution is due to the presence of the unstable phase. Therefore, if we suppose that the initial data  takes value only in the two stable phases, we can assume  that the fraction $\lambda_0\equiv 0$ and $\lambda_1$, $\lambda_2$ are proper  characteristic function. Obviously in this situation $v=\phi(u).$ 

\noindent  Starting from these considerations it was introduced the ``two--phase problem''. In this setting  we suppose that  initial data $u_0\in L^{\infty}(\R)$ satisfies

\begin{equation}\label{H2}
\left\{\begin{array}{ll}
u_0\leq b\  \mbox{in}\ (-\infty, 0),\\
u_0\geq c\  \mbox{in}\ (0, \infty),\\
\phi(u_0)\in H^1_{loc}(\R).\end{array}\right.
\end{equation}

\noindent Then, we search a solution that satisfies the entropy formulation of Plotnikov  and has a particular structure. More precisely, since the initial datum $u_0$ takes values only in the stable phases, we  impose that  solutions of the problem \eqref{B1} are again in these phases with  a regular interface separating the domain $Q_T$ into two different regions where the different phases are achieved. Finally we suppose that equation $u_t=\phi(u)_{xx}$ is fulfilled in the weak sense (see Definition \ref {two_more} in Section \ref{2}).  It is important to  underline that the entropy formulation, obtained by the pseudoparabolic  approximation \eqref{E2}, gives strong admissibility condition for the evolution of the interface that separates the two stable phases (see Proposition \ref {admis}). Existence and uniqueness results  for the two--phases problem was obtained in \cite{arma}, \cite{T}, \cite{ST}.   Then, at least for solution in the stable phases, the formulation of the forward--backward  problem suggested by Plotnikov is well--defined. In particular it is worth to note that the solutions obtained in \cite{H} do not satisfy the entropy condition introduced in \cite{Pl2}.

On the contrary, if we consider initial data in the unstable phase, it is possible to show examples of not uniqueness for solutions that satisfy the entropy formulation of Plotnikov  (see \cite{T2}). These examples are obtained by considering solutions that satisfy the forward--backward equation  in the measure--valued sense, that is  the coefficients $\lambda_i$, $i=0,1,2$ are not trivial functions. 

In this paper we want to study the two--phase problem when one of the   phases considered is the unstable one. Obviously, since we impose that the forward--backward equation is satisfied, at least in the weak sense,  we can not obtain existence for a generic initial data. However this study is interesting for different reasons.

\noindent  The first one is that this provides a  class of solution that satisfy the forward--backward problem  with unstable phase. 
Examples of explicit entropy solutions of  the forward--backward parabolic equation that takes values also in the unstable phase are given in \cite{GT}, where it  is considered the ``Riemann problem'' and   a  solution is obtained by self--similar methods.

\noindent The second one is that this study could give information about uniqueness or not uniqueness for entropy solutions of Plotnikov that still satisfy the original forward--backward equation. This is not a trivial question, since it is well--known (see e.g. \cite{Fr}) that for backward parabolic problems uniqueness results are still true. 

 This paper is organized  in two further sections.

\noindent   In Section \ref {2} we  shall state the precise definition of entropy solution.  In  particular we shall give the definition of  the two--phase problem. Moreover we shall recall which are the consequences of the entropy admissibility condition for the evolution of the interface of the two--phase stable--stable  problem. Finally we shall characterize the admissibility condition also for the two--phase problem in which one of the phase is the unstable one. 

In Section \ref{3} we shall study existence for the two--phase  stable--unstable problem.  We shall concentrate to the case in which the response function $\phi$ is piecewise linear. 

We shall prove that   explicit solutions of the two--phase  stable--unstable problem can be obtained by solving proper  generalized Abel's equations.

 \section{Entropy formulation}\label{2}

Let us recall the precise definition of entropy solution  introduced by  Plotnikov (see \cite {Pl1}, \cite{Pl2}, \cite {MTT}). Here we consider the 1-dimensional case with the domain $\Omega=\R$. 

\begin{definition}\label{defi3} 
An entropy solution  to problem (\ref{B1}) in $Q_T$
is given by  $u, \lambda_0, \lambda_1, \lambda_2 \in L^\infty (Q_T)$, $v\in L^\infty (Q_T) \cap {L^{2}((0,T), H^1_{loc}(\R))}$ such that:
\smallskip

\noindent $(a)$
$\sum\limits_{i=0}^2\lambda_i=1$, $\lambda_i\ge 0$ and there holds:
\be\label{pl1}
u=\sum_{i=0}^2\lambda_i \beta_i(v)\, 
\ee
with $ \lambda_1 =1$ if $v<A$,  $\lambda_2=1$ if $v>B$;
\smallskip

\noindent $(b)$ 
the couple $(u,v)$ is a weak solution of the  equation $u_t=v_{xx}$ in $Q_T$:
\be \label{weak}
\int \!\! \!\int_{Q_T} \!\!\!  u\psi_t  \!-\!   v \psi_x  dxdt +\int_{\Omega}u_0(x)\psi(x,0)dx=0  
\ee
for any $\psi\in C^1(\overline{Q}_T),\ \psi(\cdot,T)= 0$ in $\overline{\Omega}$.
\smallskip

\noindent $(c)$ \be \label{entrilim}
\int \!\! \!\int_{Q_T} \!\!\! \Big \{ G^*\psi_t  \!-\!  g(v) \nabla v \cdot  \nabla \psi  - g'(v)| \nabla v|^2\psi 
\Big \}\, dxdt \, \ge \, 0   \,   
\ee
for any $\psi \in C^\infty_0(Q_T)$, $ \psi \ge0$, $g \in C^1(\R)$, $g' \ge 0$ where 
\be \label{gistar}
G^*(x,t) :=\sum_{i=0}^2\lambda_iG(\beta_i(v(x,t))) \qquad \textrm{ for a.e. } (x,t) \in Q_T \, .
\ee
 and 
 \be\label{defG}
G(u):= \int^u_{k} g(\phi(s))ds, \,\,k\in\R.
\ee

 \end{definition}

We do not enter in the motivation of such definition (see e.g \cite{Pl2}, \cite{MTT}, \cite{T2}). We limit to  point out that the entropy inequality \eqref{entrilim} is satisfied by any couple of functions  $(u, v)$ obtained by proper limit of a subsequence $\{(u_{\epsilon_n},v_{\epsilon_n}) \}$, where $u_{\epsilon_n}$  is the solution of  the pseudoparabolic problem and $v_{\epsilon_n}=\phi(u_{\epsilon_n})+\epsilon_n (u_{\epsilon_n})_t$. This entropy condition implies a monotonicity condition on the  coefficients
$\lambda_i$, $i=1,2$ (see \cite{Pl2}, \cite{MTT}).

In this paper we are interested in the two--phase problem. 
We introduce the following 

\begin{definition}\label{two_more}

Let $u_0\in L^{\infty}(\R)$ be such that

$H_1$) $\phi(u_0)\in C(\R)$, $\phi(u_0)\in C^1((-\infty,0)) \cap C^1((0,\infty))$, $\phi(u_0)'\in L^1((-\infty,0)) \cap L^1((0,\infty))$;
 
 $H_2$) given the sets  $I_1=(-\infty,b]$, $I_0=[b,c]$ and $I_2=[c,+\infty)$ \begin{equation}\label{H_new}
\left\{\begin{array}{ll}
u_0\in I_i\  \mbox{in}\ (-\infty, 0),\\
u_0\in I_j\  \mbox{in}\ (0, \infty),\end{array}\right.
\end{equation}

\noindent where $i,j\in\{0,1,2\}$, $i\neq j$.

\noindent By a two--phase solution to problem \eqref{B1} we mean a triple $(u,v,\xi)$ such that:
\smallskip

\noindent $(i)$ $u\in L^{\infty}(Q_T),\, v\in L^2((0,T);H^1_{loc}(\R))$, $v\in C(\overline V_1)\cap C(\overline V_2)$, $v_x\in L^{\infty}((0,T), L^1(\R))$ and $\xi:[0,T]\to \overline{\Omega}$, $\xi\in C^1([0,T])$, $\xi(0)=0$, $v(\xi(\cdot),\cdot)\in C([0,T]))$,
$$
\lim_{\delta\to 0^+}v(\xi(t+\delta),t)=\lim_{\delta\to 0^-}v(\xi(t+\delta),t)\,\hbox{ for every  } t\in (0,T);$$
\medskip

\noindent $(ii)$ we have:
\begin{equation}\label{eq.fasi}
\left\{\begin{array}{ll}
u=\beta_i (v)\ \mbox{in}\ \,V_1,\\
u=\beta_j (v)\ \mbox{in}\ \,V_2,\end{array}\right.
\end{equation}

where
\begin{eqnarray}\label{eq.V_1}
&&V_1:=\left\{(x,t)\in \overline{Q}_T\,|\ -\infty< x<\xi(t)\,,\ t\in (0,T)\right\},\\
&&V_2:=\left\{(x,t)\in \overline{Q}_T\,|\ \xi(t)<x<+\infty\,,\ t\in (0,T)\right\},
\end{eqnarray}
and
\begin{equation}\label{eq.interface}
\gamma:=\partial V_1\cap\partial V_2=\left\{(\xi(t),t)\,|\ t\in [0,T]\right\};
\end{equation}

\noindent $(iii)$ $u$ satisfies condition $b)$  of  Definition \ref{defi3};

\noindent $(iv)$ $u$ satisfies condition  $c)$ of  Definition \ref{defi3}.

\end{definition}

By construction, a solution of the two--phase problem is also an entropy solution in the sense of Definition \ref{defi3} that has a particular structure. Observe that the previous definition implies  that  $\phi(u)=v$ everywhere so the coefficients $\lambda_k$,
 $k\in\{0,1,2\}$,  correspond to characteristic functions. More precisely
$\lambda_i=I_{V_1}$ and $\lambda_j=I_{V_2}$, where we denote with $I_E$  the characteristic function of the set $E$. In this case $G^*$ defined in \eqref{gistar} is equal to $G(u)$. 

\noindent This kind of problem was studied in different papers. However, it was evaluated  only  the case in which the initial data is in the two stable phases $I_1$ and $I_2$. In \cite {arma} and \cite{T} it was considered the piecewise linear response function $\phi$   given in \eqref{philin}. Uniqueness and local existence was obtained in \cite{arma}, global existence was stated in \cite{T}. The general nonlinear case was studied in \cite {ST}, where uniqueness was established and local existence was obtained for a class of initial data by using directly the approximation problem \eqref{E2}.

In some sense we can think that we have two different parabolic problems in the regions $V_1$ and $V_2$. These problems are related each other by \eqref{weak} and \eqref{entrilim}. The weak equation \eqref{weak} gives us the equation of the unknown interface $\xi$ and the entropy integral inequality impose some admissibility condition to the interface.  The piecewise linear case is easier to treat since in the regions $V_i$ we have classical heat equations, then we get more regularity for the solutions. In particular, we can assume the existence of the trace of the function $v_x$ along the interface. Then, using the weak equation \eqref{weak},
we obtain the following Rankine--Hugoniot equation (see  \cite{EP}, \cite{MTT})

\be\label{RH0}
\xi'(t)=\frac{v_x(\xi(t)-,t)-v_x(\xi(t)+,t)}{u(\xi(t)+,t)-u(\xi(t)-,t)},
\ee
where $v_x(\xi(t)\pm,t):=\displaystyle{\lim_{s\to\pm0}v_x(\xi(t)+ s,t)}$ and $u(\xi(t)\pm,t):=\displaystyle{\lim_{s\to\pm0}u(\xi(t)+ s,t)}.$

In the general nonlinear case it  is not evident that we can write a Rankine--Hugoniot equation for the interface. In fact  the parabolic equation in the two regions $V_1$ and $V_2$ could degenerate at values $c$ and $b$. Then  the trace of $v_x$ does not exist in the classical sense. This problem can be handled by observing that equation \eqref{weak} implies that the field $(u,-v_x)$ is divergence--free in the weak sense. Then, using the  general theory of divergence--measure vector fields, it is possible to give sense to the trace of the field along an interface. More precisely we obtain the following generalized Rankine--Hugoniot
equation
\be\label{RH}\begin{array}{ll}
&{\displaystyle\lim_{\delta\to 0^+}\frac{1}{\delta}\left[\int_{t_1}^{t_2}\!\!\!\int_{-\delta}^0\left\{\xi'(t)u(y+\xi(t),t)+v_x(y+\xi(t),t)\right\}\,\alpha(t)\,dydt\right.}\\&{\displaystyle\left.-\int_{t_1}^{t_2}\!\!\!\int^{\delta}_0\left\{\xi'(t)u(y+\xi(t),t)+v_x(y+\xi(t),t)\right\}\,\alpha(t)\,dydt\right]=0.}\end{array}
\ee
where $(t_1,t_2)\subset(0,T))$ and  $\alpha$ is any test function in  $C^1_0((t_1,t_2))$ (see \cite {ST}).

\noindent Regarding admissibility conditions for the evolution of the interface $\xi(t)$, we have to consider the entropy inequalities
\eqref {entrilim}. Let us introduce for any $g\in C^1(\R)$, $g'\ge 0$,  the   vector field $(G(u),  -[F(v)]_x)$, 
where  $G(\lambda)$ is defined in \eqref {defG} and $F(\lambda)$ is any primitive of $g(\lambda)$.  Then using \eqref{entrilim} we obtain that $(G(u),  -[F(v)]_x)$ is a measure--valued vector field. Then, again we can write a condition along the interface $\gamma$.  More precisely, for any interval $(t_1,t_2)\in (0,T)$ there holds  (see \cite{ST})
\be\label{EC}
\begin{array}{ll}&
\displaystyle{{\lim_{\delta\to 0^+}}\left(\!\frac{1}{\delta}\int_{t_1}^{t_2}\!\!\!\int_{-\delta}^0\left\{[F(v)]_x(y+\xi(t),t)+\xi'(t)\,G(u(y+\xi(t),t)\right\}\,\alpha(t)\,dydt\right.}\\&\displaystyle{-\frac{1}{\delta}\left.\int_{t_1}^{t_2}\!\!\!\int^{\delta}_0\left\{[F(v)]_x(y+\xi(t),t)+\xi'(t)\,G(u(y+\xi(t),t))\right\}\,\alpha(t)\,dydt\!\right)\le0.}\end{array}
\ee
for any test function $\alpha \in C^1_c((t_1,t_2))$, $\alpha\ge 0$.

\noindent  In some sense we can think that \eqref {EC} is an  entropy condition along the interface $\gamma$. 

Then we have less freedom to pass from one phase to the other.
This admissibility conditions are rather strong  when the two phases are the stable ones. More precisely we have (see \cite {ST} for the proof)

\begin{proposition} \label {admis} Let us consider the two--phase problem with initial data $u_0$ that satisfies \eqref{H2}.
Let $(u,v,\xi)$ be a two--phase solution of problem \eqref{B1}. Then
 \be\label{admis2}
 \xi'(t)\left\{\begin{array}{ll}
 \le 0&\hbox{ if }v(\xi(t),t)=B\\
 =0&\hbox{ if }v(\xi(t),t)\in(A,B)\\
 \ge 0&\hbox{ if }v(\xi(t),t)=A.\end{array}\right.
 \ee
 \end{proposition}

This means that the interface $\xi(t)$ moves only when $\phi(u)=v$ assumes the  critical values $A$, $B$ along it.
In particular the interface does not move when $v$ is in $(A,B)$. This allows to prove in \cite {arma}, \cite {ST} uniqueness results  for the  two--phase stable--stable problem. 

\noindent  The situation is different when one of the phase is the unstable one. We shall prove that admissibility condition is much more weaker in this situation. 

\noindent Let us consider an initial data $u_0$ that satisfies 
\be\label{st_un}\left\{\begin{array}{ll}
u_0\in (b,c)\  \mbox{in}\ (-\infty, 0),\\
u_0\in [c,\infty)\  \mbox{in}\ (0, \infty).\end{array}\right.
\end{equation}
In the next we refer to the two--phase problem with initial data  \eqref{st_un} as the ``two--phase stable--unstable problem'' and we denote it with TP S--U.
 
 Observe that we have to pay attention to the level set $v=A$. In fact,  when  $\phi$ is regular, the parabolic equations in the two regions $V_i$, $i=1,2$ have the same behavior  of the porous media equation. Then, it could exists  a region with positive measure in  which $v\equiv A$. In this case   the interface $\xi$ is not well defined since it is not clear where is the separation between the phases. For this reason in some circumstances  it could be useful to restrict the choice  of the  initial data  \eqref {st_un} such that $u_0\in (b,c)$ for  $x\le 0$. When this choice is made, it is possible to fix   $T$ small enough, such that $v(\xi(\cdot),\cdot))>A$ in $(0,T)$.

\noindent For the  TP S--U we have the following admissibility condition

\begin{proposition}\label{first_admin} Let us consider an initial data $u_0$ that satisfies \eqref{st_un}.
Then, given a  two--phase solution $(u,v,\xi)$  of the TP SU we get
$$ v(\xi(\cdot),\cdot))> A \hbox{ in }(t_1,t_2)\subset (0,T)\Longrightarrow  \xi'\le 0 \hbox{ in }(t_1,t_2)\subset (0,T).$$
\end{proposition}

\proof
We use the same techniques introduced in \cite{ST}. In order to prove the assertion by contradiction we suppose that there exists  an interval $(t_1,t_2)$ and $\rho>0$  such that $\xi'(t)>0$ and $v(\xi( t), t))>A+2\rho$  for any  $t\in (t_1,t_2)$. Then we can choose $r$ small enough such that $v\ge A+\rho$ in the domain $$S
:=\{(x,t)\in Q_T\,:\,t\in (t_1,t_2),\, \xi(t)-r\le x\le \xi(t)+r\}.$$ Let us choice an increasing  function $g$
such that
$$g=\left\{\begin{array} {ll} <0&\hbox{ in } (-\infty, A+\rho)\\ 0&\hbox{ in } [A+\rho, \infty).\end{array}\right.$$

\noindent Moreover we consider
$$G(s)=\int_{u^{\rho}}^s g(\phi(\tau)\,d\tau,\,\,\,\, F(s)=\int_B^s g(\tau)\,d\tau,$$
where $u^{\rho}$ is the point in $(c,+\infty)$ such that $\phi(u^{\rho})=A+\rho$.

Then by the entropy condition \eqref{EC} we obtain 
\be\label{entropie1}\begin{array}{ll}&
{\displaystyle\int_{t_1}^{t_2}\{G(u_1(t))-G(u_2(t))\}\,\xi'(t)\alpha(t)\,dt}\\
&+{\displaystyle\lim_{\delta\to0^+}\left(\frac{1}{\delta}\int_{t_1}^{t_2}\left\{F(v(\xi(t),t))-F(v(\xi(t)-\delta,t))\right\}\,\alpha(t)\,dt\right.}\\
&-{\displaystyle\left.\frac{1}{\delta}\int_{t_1}^{t_2}\left\{F(v(\xi(t)+\delta,t))-F(v(\xi(t),t))\right\}\,\alpha(t)\,dt\right)\le 0.}
\end{array}
\ee
Where $u_1(t):=\displaystyle{\lim_{s\to 0^-}u(\xi(t)+s,t)}$ and $u_2(t):=\displaystyle{\lim_{s\to 0^-}u(\xi(t)+s,t)}$.
Observe that by construction we have $F(v(\xi(t)\pm\delta,t))=F(v(\xi(t),t))=G(u_2(t))=0$ for $\delta\le r$.  Moreover,   for every $t\in (t_1,t_2)$,  the function $g(\phi(\cdot))$ assume negative value in a proper subinterval of  $(u_1(t),u^{\rho})$ , then we have $G(u_1(t))=\int _{u^{\rho}}^{u_1(t)}g(\phi(s))\,ds>0$. Therefore, using $\xi'(t) >0$,  inequality  \eqref{entropie1} gives a contradiction.
\qed

In the next Section we shall focus to the solution of the TP S--U in the case of the piecewise linear function $\phi$ given in \eqref{philin}. As observed in the introduction, in this context we may assume more regularity for the solution of the two phase problem. More precisely, in the piecewise linear case,  we generalize the  definition of the two--phase ``stable--stable'' problem given in \cite{arma}, \cite{T}. 

\begin{definition}\label{two_lin}  Let  $\phi$ be the piecewise linear function \eqref{philin}. Let $u_0$ be an initial data that satisfies the hypothesis  $H_1)$, $H_2)$  of Definition \ref  {two_more}.  By a two--phase regular solution of problem \eqref {B1} we mean a triple $(u,v,\xi)$ that is solution in the sense of Definition \ref{two_more} that satisfy also

$I)$\, $\xi\in C^{\frac32}(0,T);$

$II)\,$ $u \in C^{2,1}\big(V^1\big)\cap C^{2,1} \big (V_2 \big)$,
$u_x \in L^{\infty}(Q_{T})$, and  for any $t \in (0,T)$ there exist finite the limits
\begin{equation*}
\lim_{\eta \to 0^+} u_x(\xi(t) \pm \eta,t)=: u_x(\xi(t)^\pm,t).
\end{equation*}
\end{definition}
This definition is motivated in the piecewise linear case  by the existence results obtained in \cite{arma} for the two--phase stable--stable problem.

Then, we can prove the following admissibility result for the TP S--U

\begin{proposition} \label{xineg}Let  $\phi$ be the piecewise linear function \eqref{philin}. Let $u_0$ satisfies the hypothesis of Definition \ref{two_more} with $u_0\in I_0$ for $x< 0$ and $u_0\in I_2$ for $x\ge 0$. Then for any two--phase regular solution  $(u,v,\xi)$ of the TP S--U such that $u\not\equiv c$ in $Q_T$ we have $\xi'\le 0$ in $(0,T)$.   

\end{proposition}  

\proof
Let us suppose that there exists $\overline t$ such that $\xi'(\overline t)>0$. Then if $v(\xi(\overline t),\overline t)>A$ we obtain a contradiction by using Proposition  \ref{first_admin}. Therefore $v(\xi(\overline t),\overline t)=A$. Let us assume that there exist a sequence $t_n$ converging to $\overline t$ such that $v(\xi(t_n),t_n)>A$. Again by continuity we have a contradiction. Then the only possibility is that there exists an interval $(t_1,t_2)$ containing $\overline t$ such that $\xi'(\cdot)>0$ and $u(\xi(\cdot),\cdot)\equiv c$ in $(t_1,t_2)$. Suppose that there exists $s_2\in (t_1,t_2)$ such that $u_x(\xi(s_2)^+,s_2)=0$. Then using the well--known results on the maximum principle  (see e.g \cite {Fr}), we get that there exist $\overline s_2\in (t_1,t_2)$ such that  $u\equiv c$ in $V_2\cap \R\times [0,\overline s_2]$. Analogously, if there exists $s_1\in (t_1,t_2)$ such that $u_x(\xi(s_1)^+,s_1)=0$, then there exist $\overline s_1\in (t_1,t_2)$ such that  $u\equiv c$ in $V_1\cap \R\times [\overline s_1,T]$. Therefore we deduce that, if $u\not\equiv c$ in $Q_T$, necessarily $ u_x(\xi(t)^+,t)>0$ or   $u_x(\xi(t)^-,t)>0$ for every  $t\in (t_1,t_2)$.

Let us suppose that $u\not\equiv c$ in $Q_T$.  This  contradict the Rankine--Hugoniot condition \eqref{RH}. In fact we get

\begin{eqnarray*}\label{RH2}
&&0=\int_{t_1}^{t_2}\alpha(t)\xi'(t)(c-c)\,dt\\
&&=\lim_{\delta\to 0^+}\frac{1}{\delta}\left (\int_{t_1}^{t_2}\alpha(t)((v(\xi(t)-\delta,t)-A)+(v(\xi(t)+\delta,t)-A))\,dt\right)\\
&&\int_{t_1}^{t_2}\alpha(t)((\gamma_2u_x(\xi(t)^+,t)-\gamma_0u_x(\xi(t)^-,t))\,dt
\end{eqnarray*}
 for any nonnegative test function $\alpha(t)\in C^1_c((t_1,t_2))$.  
 
 \noindent Since $\gamma_0<0$ and $\gamma_2>0$  (the coefficient given in \eqref{philin}) the left hand side in \eqref{RH2} is strictly positive then we get the contradiction.

\qed 

In the following we always assume that $\phi$ is the piecewise function defined in \eqref{philin}.

Using  Proposition \ref{xineg},  we see that the admissibility condition for the TP S--U  implies  that the unstable phase tends to disappear respect to the stable one. In the next Proposition we see that the condition $\xi'\le 0$ is in fact equivalent  to the admissibility condition.
\begin{proposition}\label{xiifonlyif} Let $u_0$ satisfies the hypothesis  of Proposition \ref{xineg}. Let  $(u,v,\xi)$ be a triple that satisfies conditions $i)$, $ii)$, $iii)$ of Definition \ref{two_more} and condition $I)$ and $II)$ of Definition \ref{two_lin}, moreover assume that $u\not\equiv c$ in $Q_T$. Then the triple   $(u,v,\xi)$ is a regular solution of the TP S--U  if and only if $\xi'\le 0$ in $(0,T)$.
\end{proposition} 

\proof Using Proposition \ref {xineg} it remains to prove that 
$\xi'\le 0$ in $(0,T)$ implies that the triple $(u,v,\xi)$ satisfies the entropy admissibility condition \eqref{entrilim}.

\noindent First of all, we observe that, by  assumption, the solution is regular in the regions $V_1$ and $V_2$. Therefore the entropy admissibility condition \eqref{entrilim} is equivalent to  the admissibility condition \eqref{EC} along the boundary.
 Moreover, there exists the trace of the functions $u$ and $v$ along the interface $\gamma$. Let us  introduce  the notations
 $u_1(t)=u(\xi(t)^-,t)$, $u_2(t)=u(\xi(t)^+,t)$, $v_x^1(t)=v_x(\xi(t)^-,t)$, $v_x^2(t)=v_x(\xi(t)^+,t)$. Therefore condition \eqref{EC} becomes
\be\label{entropie2}\begin{array}{ll}&
{\displaystyle\int_{0}^{T}\{G(u_1(t))-G(u_2(t))\}\,\xi'(t)\alpha(t)\,dt}\\
&+{\displaystyle\int_{0}^{T}\left(g(\phi(u_1(t)))v_x^1(t)-g(\phi(u_2(t)))v_x^2(t)\right)\alpha(t)\,dt\le 0.}
\end{array}
\ee
for any test function $\alpha \in C^1_c((0,T))$, $\alpha\ge 0$. Let us recall that, for hypothesis, $\phi(u_1(t))=\phi(u_2(t))$ for any $t\in (0,T)$.  

\noindent Then, in order to prove \eqref{entropie2} it is enough to check  that 
\be\label{h_t} h(t):=(G(u_1(t))-G(u_2(t)))\xi'(t)+\phi(u_2(t))(v_x^1(t)-v_x^2(t))\le 0\ee 
for any $t\in (0,T)$. 

\noindent Let us observe that the Rankine--Hugoniot condition \eqref{RH} becomes in this ``regular'' case 

\be\label{RH3}
\int_{0}^{T}\beta(t)\xi'(t)(u_1(t)-u_2(t))\,dt-
\int_{0}^{T}\beta(t)((v_x^2(t)-v_x^1(t))\,dt=0
\ee
 for any test function $\beta \in C^1_c((0,T))$. Let us consider the set $J:=\{t\in (0,T)\,:\, u_1(t)\neq u_2(t)\}$. We deduce 
  the classical Rankine--Hugoniot condition
 \be\label{RH4}
 \xi'(t)=\frac{v_x^2(t)-v_x^1(t)}{u_1(t)-u_2(t)}
 \ee
for any $t\in J$. 

\noindent  Therefore, using $\xi'(t)\le 0$ and $\phi(s)\le \phi(u_2(t))=\phi(u_1(t))$ for any $s\in (u_2(t))-u_1(t)))$,   we get

$$
h(t)=-\xi'(t)\left(G(u_2(t))-G(u_1(t))-g(\phi(u_2(t)))(u_2(t)-u_1(t))\right)=
$$
$$
-\xi'(t)\left(\int_{u_1(t)}^{u_2(t)}g(\phi(s))-g(\phi(u_2(t)))\,ds\right)\le0
$$
for any $t\in J$. Let us prove that $\overline J=[0,T]$. If this is not true there exists an interval $I=(t_1,t_2)$ such that $u_1\equiv u_2\equiv c$ in  $I$. Then, reasoning as in the proof of Proposition \ref{xineg}, we obtain $u\equiv c$ in $Q_T$ that contradicts the hypothesis.

\qed

 \section{Existence results for  two phase stable--unstable problem}\label{3}
 In this Section we analyze the existence of solutions of the TP S--U. Obviously we can not expect to have solution for general initial data. In fact, in the left hand side of the interface $\gamma$, we have to consider  a backward parabolic problem. In particular  $u_0$ for $x\le 0$ has to be  much more regular respect to the hypothesis $H_1)$  in Definition \ref{two_more}. 
 
\noindent  Here we are interested in obtaining  explicit solutions of the TP S--U without imposing the initial condition in the semiaxis  $x\le 0$. This problem is not well-posed since uniqueness could not be true. In any case it is interesting  to study the  inverse problem in which we impose the final data $u(\cdot,T)$ in the backward region.

Therefore we search regular solutions of the two phase stable--unstable problem that satisfy all the conditions of Definition \ref{two_lin} but we replace initial condition   $u(x,0)=u_0(x)$ for $x\le 0$ with $u(x,T)=u_T(x)$ for $x\le K$, where $K$ is a not positive constant and $u_T$ is a given proper regular function defined in $(-\infty,L)$ that takes values in the unstable phase. 

\noindent In the following we assume more regularity for the interface $\xi$ and data $u_0$ and $u_T$.
More precisely we give the following definition 

\begin{definition}\label{u_T} Let $K$ be a fixed constant in $(-\infty,0]$. Let $u_0\in L^{\infty}(\R^+)$ such that 
$u_0\in C^2(\R^+)$, $u_0'\in L^1(\R^+)$, $u_0''\in L^{\infty}(\R^+)$, $u_0\ge c$ in $\R$. Let $u_T\in C^2((-\infty,K))$, $u_T'\in L^1(-\infty,K)$, $u_T''\in L^{\infty}(-\infty,K)$ $u_T\in (b,c)$ in $(-\infty,K]$. 
By a regular two phase entropy solution of the problem
 \be\label{B_T}
 \left\{\begin{array}{lll} u_t= \phi(u)_{xx} & \hbox{ in } Q_T=\R\times (0,T)\\
  u(x,0)=u_0(x) & \hbox{ in }\R^+\times\{0\}\\
  u(x,T)=u_T(x) & \hbox{ in }(-\infty,K)\times\{0\},\end{array}
  \right.
 \ee

 we mean a triple $(u,v,\xi)$ such that:
\smallskip

\noindent $(i)$ $u\in L^{\infty}(Q_T),\, v\in L^2((0,T);H^1_{loc}(\R))$, $v\in C(\overline V_1)\cap C(\overline V_2)$, $v_x\in L^{\infty}((0,T), L^1(\R)))$ and  $\xi\in C^{1+\beta}([0,T])$, $\beta>\frac12$, $\xi(0)=0$, $\xi(T)=K$, $v(\xi(\cdot),\cdot)\in C([0,T]))$,
$$
\lim_{\delta\to 0^+}v(\xi(t+\delta),t)=\lim_{\delta\to 0^-}v(\xi(t+\delta),t)\,\hbox{ for every  } t\in (0,T);$$
\medskip

\noindent $(ii)$ we have:
\begin{equation}\label{eq.fasi2}
\left\{\begin{array}{ll}
u=\beta_0 (v)\ \mbox{in}\ \,V_1,\\
u=\beta_2 (v)\ \mbox{in}\ \,V_2,\end{array}\right.
\end{equation}

where
\begin{eqnarray}\label{eq.V_1_2}
&&V_1:=\left\{(x,t)\in \overline{Q}_T\,|\ -\infty< x<\xi(t)\,,\ t\in (0,T)\right\},\\
&&V_2:=\left\{(x,t)\in \overline{Q}_T\,|\ \xi(t)<x<+\infty\,,\ t\in (0,T)\right\},
\end{eqnarray}
and
\begin{equation}\label{eq.interface2}
\gamma:=\partial V_1\cap\partial V_2=\left\{(\xi(t),t)\,|\ t\in [0,T]\right\};
\end{equation}

 $u \in C^{2,1}\big(V^1\big)\cap C^{2,1} \big (V_2 \big)$,
$u_x \in L^{\infty}(Q_{T})$, and  for any $t \in (0,T)$ there exist finite the limits
\begin{equation*}
\lim_{\eta \to 0^+} u_x(\xi(t) \pm \eta,t)=: u_x(\xi(t)^\pm,t).
\end{equation*}

\noindent $(iii)$ the function  $u$ fulfills   the equation $u_t=\phi(u)_{xx}$ in the regions $V_0$ and $V_2$;
with $u(\cdot,0)=u_0$ in $\R$ and $u(\cdot,T)=u_T$ in $(-\infty, K)$;

\noindent $(iv)$ the   Rankine--Hugoniot condition \eqref{RH2} is satisfied;

\noindent $(iv)$ $\xi'(t)\le 0$ for every $t\in (0,T)$.

\end{definition}

\begin{remark} Using Proposition \ref{xiifonlyif} it follows that a two phase entropy solution of problem \eqref{B_T} it is  a two phase entropy solution of problem \eqref{B1} with initial data
$$
\overline u_0=\left\{\begin{array}{lll}u_0&\hbox{ in } \R^+\\
u(\cdot,0)\,&\hbox{ in } \R^-.\end{array}\right.   
$$
\end{remark}
\begin{remark} If we assume that $v(\xi(\cdot),\cdot)\neq A $ in $(0,T)$, we can replace the integral Rankine--Hugoniot condition  with the more classical  condition
\eqref{RH0}.\end{remark}

 The introduction of the previous problem could be useful when considering inverse problem. However in this paper the unique scope is to obtain explicit entropy solution
 of the forward--backward parabolic equation in the case in which there is   the unstable phase  with one of the stable ones.
 
  We shall see, in general, that if a solution of  the two phase problem \eqref{B_T} exists  it is not unique. In order to study this problem we fix a decreasing function $\xi\in C^{1+\beta}((0,T))$, $\beta>\frac 12$ and we search an entropy two phase solution of problem  \eqref{B_T} with $K=\xi(T)$. Moreover we suppose   that 
  $$\displaystyle{\sup_{x\in(-\infty,K)} u_T(x)<c},\,\displaystyle{\inf_{x\in(0,\infty)} u_0(x)>c}$$ and we search a solution such that $v(\xi(\cdot),\cdot)\neq A $ in $(0,T)$.
  
  In order to study the two phase problem we need some results on the Dirichlet--to--Neumann map in Time--dependent Domains proved  in  \cite {FP} (see also \cite{Fo}). More precisely 
  
  \begin{theorem}\label{fokas}  Let $\xi(t)\in C^1((0,T))$, $q_0\in C^1([0,\infty))$ such that $q_0'\in L^1((0,\infty))$, $g_0\in C^1([0,T])$. Let us denote with $q(x,t)$ the solution of the Dirichlet problem
  \be\left\{\begin{array}{ll}
  q_t=q_{xx} &\hbox{ in } V_{\xi}\\
  q(\xi(t),t)=g_0(t)&\hbox{ in } (0,T)\\
  q(x,0)=q_0(x)&\hbox{ in } \R
  \end{array}\right.\ee
  where $V_{\xi}:=\left\{(x,t)\in \overline{Q}_T\,|\ \xi(t)<x<+\infty\,,\ t\in (0,T)\right\}$. Then the function $f(t):=q_x(\xi(t),t)$ is characterized as the solution of the following Volterra equation
  \be
  \pi f(t)=N(t)+\int_0^t K(s,t)f(s)\,ds
  \ee
  where the integral kernel is defined by
  $$
  K(s,t)=\frac{\sqrt{\pi}}{2}\frac{\xi(t)-\xi(s)}{t-s}
  \frac{e^{-\frac{(\xi(t)-\xi(s))^2}{4(t-s)}}}{(t-s)^{\frac12}}\,\, 0<s<t<T
  $$
  and the function $N(t)$ that depends on the initial--boundary value  is given by
 $$ N(t)=\sqrt{\pi}\left[\frac{1}{\sqrt{t}}\int_0^{\infty}e^{-\frac{(\xi(t)-x)^2}{4t}}q_0'(x)\,dx-\int_0^t  \frac{e^{-\frac{(\xi(t)-\xi(s))^2}{4(t-s)}}}{(t-s)^{\frac12}}g_0'(s)\,ds\right]
 $$
  \end{theorem}

We shall obtain a solution of the two phase problem \eqref{B_T} by imposing the Rankine--Hugoniot condition \eqref{RH0},  the continuity of the solution $v=\phi(u)$ along the interface. Then by the characterization obtained in Theorem   \ref{fokas} we shall see that the solution of the original problem is associated to a generalized Abel problem. 

We proceed as follows. We have to consider the following parabolic problems

  \be\label{P_1}\left\{\begin{array}{ll}
  u_t=\gamma_0u_{xx} &\hbox{ in } V_1\\
  q(\xi(t),t)=g_0^-(t)&\hbox{ in } (0,T)\\
  q(x,T)=u_T&\hbox{ in } (-\infty,K)
  \end{array}\right.\ee

 \be\label{P_2}\left\{\begin{array}{ll}
  u_t=\gamma_2u_{xx} &\hbox{ in } V_2\\
  q(\xi(t),t)=g_0^+(t)&\hbox{ in } (0,T)\\
  q(x,0)=u_0&\hbox{ in } \R
  \end{array}\right.\ee
where $V_1$ and $V_2$ are given in \eqref{eq.V_1}, $g_0^+$, $g_0^-$ are unknown of the problem that are related each other by the continuity condition $\phi(u(\xi(t)-,t))=\phi(u(\xi(t)+,t))$ in $(0,T)$. More precisely, we have
$$ \gamma_0 g_0^-(t)+\delta_0=\gamma_2 g_0^+(t)+\delta_2=:m(t).$$

\noindent Finally we determinate the unknown function $m(t)$ by using the  Dirichlet--Neumann maps associated to the problems \eqref{P_1} and \eqref{P_2} and imposing the Rankine--Hugoniot condition \eqref{RH0}. We shall obtain such  Dirichlet--Neumann maps by standard change of variable. More precisely we have
\begin{proposition} \label{DN1_map}Let $u^+$ be the solution of the parabolic problem \eqref{P_2}. Then the function $f^+(t)=u^+_x(\xi(t),t)=f_2(\gamma_2 t)$, where $f_2$ is the solution of the following Volterra problem 
\be\label{volterra2}
  \pi f_2(\tau)=N_2(\tau)+\int_0^{\tau} K_2(s,\tau)f_2(\tau)\,ds,\,\,\,\, 0\le\tau\le\gamma_2 T
  \ee
with the kernel
$$
K_2(s,\tau)=\frac{\sqrt{\pi}}{2}\frac{\overline\xi(\tau)-\overline \xi(s)}{\tau-s}
  \frac{e^{-\frac{(\overline \xi(\tau)-\overline\xi(s))^2}{4(\tau-s)}}}{(\tau-s)^{\frac12}}\,\, 0<s<\tau<\gamma_2T,
  $$
$\overline \xi(\tau)=\xi(\frac{\tau}{\gamma_2})$ and 
$$N_2(\tau)=\sqrt{\pi}\left[\frac{1}{\sqrt{\tau}}\int_0^{\infty}e^{-\frac{(\overline \xi(\tau)-x)^2}{4\tau}}u_0'(x)\,dx-\int_0^\tau \frac{e^{-\frac{(\overline \xi(\tau)-\overline \xi(s))^2}{4(\tau-s)}}}{(\tau-s)^{\frac12}}\frac{{g_0^+}'(\frac{s}{\gamma_2})}{\gamma_2}\,ds\right]
 $$
 \end{proposition}
\proof This is consequence of Theorem \ref{fokas}. 

\noindent Let us  consider the function $w(y,\tau):=u^+(y,\frac{\tau}{\gamma_2})$ then this is solution of the problem 

 \be\label{P_2bis}\left\{\begin{array}{ll}
  w_{\tau}=w_{yy} &\hbox{ in } \overline V_2\\
  w(\overline \xi(\tau),\tau)=g_0^+(\frac{\tau}{\gamma_2})&\hbox{ in } (0,\gamma_2T)\\
  w(y,0)=u_0(y)&\hbox{ in } \R
  \end{array}\right.\ee
   where $\overline \xi(\tau)=\xi(\frac{\tau}{\gamma_2})$ is defined in $(0,\gamma_2 T)$ and $$\overline V_2=\left\{(x,\tau)\in \overline{Q}_{\gamma_2 T}\,|\ \overline \xi(\tau)<x<+\infty\,,\ \tau\in (0,\gamma_2 T)\right\}.$$
   
   \noindent Therefore $f_2(\tau)=w_x(\overline \xi(\tau),\tau)$ is characterized by Theorem \ref{fokas}. Obviously $f^+(t)=f_2(\gamma_2 t)$ then we obtain the thesis.
   \qed
   
   Analogously we can obtain the Dirichlet--Neumann map for the parabolic problem \eqref {P_1}. More precisely 
   
   \begin {proposition}\label{DN2_map} Let $u^-$ be the solution of the parabolic problem \eqref{P_1}. Then the function $f^-(t)=u^-_x(\xi(t),t)=-f_1(|\gamma_0|(T- t))$, where $f_1$ is the solution of the following Volterra problem 
\be\label{volterra1}
  \pi f_1(\tau)=N_1(\tau)+\int_0^{\tau} K_1(s,\tau)f_1(\tau)\,ds,\,\,\,\, 0\le\tau\le|\gamma_0| T
  \ee
with the kernel
$$
K_1(s,\tau)=\frac{\sqrt{\pi}}{2}\frac{\underline\xi(\tau)-\underline \xi(s)}{\tau-s}
  \frac{e^{-\frac{(\underline \xi(\tau)-\underline\xi(s))^2}{4(\tau-s)}}}{(\tau-s)^{\frac12}}\,\, 0<s<\tau<|\gamma_0|T,
  $$
$\underline  \xi(\tau)=K-\xi(T-\frac{\tau}{|\gamma_0|})$ and 
$$N_1(\tau)=\sqrt{\pi}\left[-\frac{1}{\sqrt{\tau}}\int_0^{\infty}e^{-\frac{(\underline \xi(\tau)-x)^2}{4\tau}}u_T'(K-x)\,dx+\int_0^\tau \frac{e^{-\frac{(\underline \xi(\tau)-\underline \xi(s))^2}{4(\tau-s)}}}{(\tau-s)^{\frac12}}\frac{{g_0^-}'(T-\frac{s}{|\gamma_0|})}{|\gamma_0|}\,ds\right].
 $$
\end{proposition}
   
   \proof In this case we consider the function $\theta(y,\tau):=u^-(K-y,T-\frac{\tau}{|\gamma_0|})$. It is easy to check that $\theta$ satisfies the following parabolic problem
   
 \be\label{P_1bis}\left\{\begin{array}{ll}
  \theta_{\tau}=\theta_{yy} &\hbox{ in } \overline V_1\\
  \theta(\underline \xi(\tau),\tau)=g^-_0(T-\frac{\tau}{|\gamma_0|})&\hbox{ in } (0,|\gamma_0|T)\\
  \theta(y,0)=u_0(K-y)&\hbox{ in } \R
  \end{array}\right.\ee
where    $$\overline V_1=\left\{(x,\tau)\in \overline{Q}_{|\gamma_0| T}\,|\ \underline \xi(\tau)<x<+\infty\,,\ \tau\in (0,|\gamma_0| T)\right\}.$$ 
Again we use Theorem \ref{fokas} and the relation $f^-(t)=-f_1(|\gamma_0|(T-t))$ with $f_1(\tau)=\theta_x(\underline\xi(\tau),\tau)$ to obtain the thesis.
   
   \qed
\begin{remark} In order to determinate the unknown $m(t)=v(\xi(t),t))$ we impose the Rankine--Hugoniot condition \eqref{RH0} that for Propostitions \ref{DN1_map}, \ref{DN2_map} becomes
\be\label{RH_new}
\xi'(t)=\frac{|\gamma_0|f_1(|\gamma_0|(T-t))-\gamma_2 f_2(\gamma_2 t)}{g^+(t)-g^-(t)}.
\ee 
We shall prove that equation \eqref{RH_new} could be written as a generalized Abel equation.
\end{remark}

First of all we have to enunciate this  result on the Volterra equation (for a proof see e.g. \cite{T}).
\begin{proposition}\label{prop_volterra}
\label{volt_rad}
Let   $K$ be a continuous function defined in $C_T:=\{(t,s)\in (0,T)\times(0,T)\,:\, 0<s<t\}$ and $M$ be a positive constant such that  $$|K (t,s)|\le \frac{M}{\sqrt{t-s}} \qquad\qquad\hbox{for every } (t,s)\in C_T.$$

Let us consider the operator $L:C([0,T])\rightarrow C([0,T])$ defined by

\be\label{volterra}
L(x(t))=x(t)-\int_0^t K(t,s)x(s)\,ds.
\ee
Then, there exist a continuous function $\mathcal{H}(t,s)$ defined in $C_T$ and a constant $\overline M>0$,
such that 
\be\label{abel}
x(t)=L(x(t))+\int_0^t \mathcal{H}(t,s)L(x(s))\,ds,
\ee
with  $\mathcal{H}(t,s)=K(t,s)+G(t,s)$ in $C_T$
and $G$ continuous in $\overline C_T$.\end{proposition}
\begin{remark}\label{G_alfa} Following the proof given in \cite{T}  (Lemma 2.7) we obtain that 
\be\label{G_def}
G(t,s)=\sum_{i=2}^{+\infty}H_i(t,s)
\ee
where
\be\label{k1}
H_1(t,s)=K(t,s)I_{(0,t)}(s),
\ee
 we denote, again, with $I_E$ the characteristic function of the set $E$
and $H_n$, $n\ge 2$ are defined recursively by
$$H_n(t,s)=\int_0^T H_{n-1}(t,z)H_1(z,s)\,dz.$$
Moreover the convergence in \eqref{G_def} is uniform in $[0,T]\times[0,T]$.

\noindent More precisely  in \cite {T}  it is proved by induction that
\be\label{indu}
|H_{n+2}(t,s)|\le \frac{M^{n+2}\pi^{\frac {n+1}2}\Gamma(\frac 12)}{\Gamma(\frac{2+n}{2})}(t-s)^{\frac{n}{2}}I_{(0,t)}(s)
\ee
 for every $n\in\N$, 
where $\Gamma$ is the classical Gamma function.  

\noindent Let us put $ h_1(t,s)=H_1(t,s)\sqrt{t-s}$,  suppose that there exist $\alpha<\frac 12$ and $S>0$ such   that 

\be\label{h_1}
|\partial_t h_1(t,s)|, |\partial_s h_1(t,s)|\le \frac S{{(t-s)^{\alpha}}}\hbox{ in }C_T
\ee  then  using the same technics that give the estimates \eqref{indu} we obtain that there exist a constant $\overline S$ such that 
\be\label{GGG}|\partial_t G(t,s)|, |\partial_s G(t,s)|\le \frac {\overline S}{{(t-s)^{\frac12}}}\hbox{ in }C_T.\ee
\end{remark}

Using Proposition \ref{prop_volterra} we can give the functions  $f_1$ and $f_2$ in the Volterra equation \eqref{volterra2} and \eqref{volterra1} in a more explicit way. More precisely we have

\begin{proposition} \label{prop_kilbas} There exists an Holder continuous  function   $G_2(t,s)$, of order $\frac 12$, defined in $\overline C_{\gamma_2 T}$, satisfying \eqref {GGG},   such that 
\be
f_2(t)=\frac{N_2(t)}{\pi}+\frac{1}{\pi^2}\int_0^t (K_2(t,s)+G_2(t,s)) N_2(s)\,ds\,\hbox{ in } (0,\gamma_2 T)
\ee
where $f_2$ and $N_2$ are the functions defined in Proposition \ref{DN2_map}. 

\end{proposition}

\proof It is enough to use Proposition \ref{prop_volterra} and Remark \ref {G_alfa}. In fact 
 by  a straightforward  calculation we obtain that there exists a constant  $M_2$ such that
$$
K_2(s,t)\le \frac{M_2}{\sqrt{t-s}}\hbox{ in }C_{|\gamma_2|T}.$$ Moreover since  $\xi\in C^{1+\beta}$, $\beta>\frac12$ we obtain that $h_1(t,s)=\sqrt{t-s}K_2(t,s)$ satisfies estimate \eqref{h_1} with $\alpha=1-\beta<\frac12$. Therefore the thesis is consequence of the equation \eqref{volterra2}
\qed

\noindent Analogously we get

\begin{proposition} There exists an Holder continuous  function   $G_1(t,s)$, of order $\frac 12$, defined in  $\overline C_{|\gamma_0| T}$, satisfying \eqref {GGG}, such that 
\be
f_1(t)=\frac{N_1(t)}{\pi}+\frac{1}{\pi^2}\int_0^t (K_1(t,s)+G_1(t,s)) N_1(s)\,ds\,\hbox{ in } (0,|\gamma_0| T)
\ee
where $f_1$ and $N_1$ are the functions defined in Proposition \ref{DN1_map}.

\end{proposition}

Using the previous results  the Rankine--Hugoniot equation \eqref{RH_new} becomes
\be\label{afraid}\begin{array}{ll}
&\displaystyle{\xi'(t)\left(\frac{m(t)-\delta_2}{\gamma_2}-\frac{m(t)-\delta_0}{\gamma_0}\right)=}\\
&+\displaystyle{\frac{|\gamma_0| N_1(|\gamma_0|(T- t))}{\pi}+\frac{|\gamma_0|}{\pi^2}\!\!\int_0^{|\gamma_0|(T- t)} \!\!\!\!\!\!\!\!\!\!\!\!\!\!\!\!(K_1(|\gamma_0|(T- t),s)+G_1(|\gamma_0|(T- t),s))N_1(s)\,ds}\\&-\displaystyle{\frac{\gamma_2N_2(\gamma_2 t)}{\pi}-\frac{\gamma_2}{\pi^2}\int_0^{\gamma_2 t}\!\!\!\! (K_2(\gamma_2 t,s)+G_2(\gamma_2 t,s)) N_2(s)\,ds}\end{array}
\ee
for every $t\in (0,T)$. Let us observe that $N_1$ and $N_2$ depends on the unknown $m(t)$. Then we have to analyze in detail the previous equation in order to state that this is equivalent to an Abel's type equation for the function $m'(t)$. This is the content of the following

\begin{theorem} Equation \eqref{afraid} is equivalent to the following generalized
  Abel's equation of the first order:
\be\label{abel1}
\int_0^t\frac{k_1(t,s)}{\sqrt{t-s}}{m'(s)}\,ds+\int_t^T\frac{k_2(t,s)}{\sqrt{s-t}}{m'(s)}\,ds=h(t)
\ee
where $k_1$, respectively $k_2$,  is an Holder continuos function of order $\frac 12$  defined in the set
$\overline C_T$, respectively $[0,T]\times[0,T]\setminus C_T$ and $h(t)$ is a function in $C^{\frac12}([0,T])$  that depends on the initial data $u_0'$, $u_T'$ and $m(0)$.
\end{theorem}

In order to prove the previous Theorem we introduce the following lemmas.

\begin{lemma}\label{lemma_kilbas0} Let $g_1$, $g_2$ Holder continuous function of order $\gamma$ defined in $C_T$.
Let us consider the function
$$
h(t,s)=\int_s^t\frac{g_1(t,z)g_2(z,s)}{\sqrt{t-z}\sqrt{z-s}}\,dz
$$
then $h$ is Holder continuous of order $\gamma$ in $C_T$.
\end{lemma}
Let us fix $(t,s)\in C_T$,  consider for every $\tau\in (s-t,T-t)$ the following   change of variable $w=f(w)=s+(w-s)\frac{t-s+\tau}{t-s}$. Then we get
\be\label{g_1g_2}
h(t+\tau,s)=\int_s^t\frac{g_1(t+\tau,f(w))g_2(f(w),s)}{\sqrt{t-w}\sqrt{w-s}}\,dw.
\ee
Therefore
$$h(t+\tau,s)-h(t,s)=\int_s^t\frac{g_1(t+\tau,f(w))g_2(f(w),s)-g_1(t,w)g_2(w,s)}{\sqrt{t-w}\sqrt{w-s}}\,dw.$$
Let us estimate the numerator in the integrand of \eqref{g_1g_2}. We get

$$|g_1(t+\tau,f(w))g_2(f(w),s)-g_1(t,w)g_2(w,s)|\le$$
$$ |g_2(f(w),s)||g_1(t+\tau,f(w))-g_1(t,w)|+|g_1(t,w)||g_2(f(w),s)-g_2(w,s)|\le$$
$$
|g_2(f(w),s)|(|g_1(t+\tau,f(w))-g_1(t,f(w))|+|g_1(t,f(w))-g_1(t,w)|)+$$
$$
|g_1(t,w)||g_2(f(w),s)-g_2(w,s)|.
$$
Using the Holder properties of $g_1$ and $g_2$ we obtain 
$$
|g_1(t+\tau,f(w))g_2(f(w),s)-g_1(t,w)g_2(w,s)|\le P (|\tau|^{\gamma}+|f(w)-w|^{\gamma})=$$
$$
P \left (|\tau|^{\gamma}+\left|\frac {\tau (w-s)}{t-s}\right|^{\gamma}\right)\le 2P|\tau|^{\gamma},
$$
where $P$ is a constant that does not depend on $s,t,\tau$. This implies that $h(\cdot,s)$ is Holder continuous of order $\gamma$ uniformly in $s$. Analogously we can prove that $h(t,\cdot)$ is Holder continuous of order $\gamma$ uniformly in $t$.\qed

\begin{lemma}\label{lemma_kilbas} Let $\xi\in C^{1+\beta}([0,T])$, $\beta>\frac{1}{2}$, $f\in C^{\frac12}([0,T])$. Let us consider the following function
$$
g(t):=\int_0^t K(t,s)f(s)\,ds,
$$
where
\be\label{KeG}
K(t,s)=\frac{\xi(t)- \xi(s)}{t-s}
  \frac{e^{-\frac{( \xi(t)-\xi(s))^2}{4(t-s)}}}{(t-s)^{\frac12}}.
\ee
Then $g\in C^{\frac 12}([0,T]).$
\end{lemma}

\proof Let  $t_1, t_2\in [0,T]$. It is not restrictive to assume $t_1< t_2$. By changing variable, we can write 
\be\label{g_lip}
g(t)=\int_0^t \frac{\xi(t)-\xi(t-y)}{y^{\frac 32}} e^{-\frac{( \xi(t)-\xi(t-y))^2}{4y}}f(t-y)\,dy.
\ee
In order to simplify the notations we introduce the functions
$$g_1(t,y)=\frac{\xi(t)-\xi(t-y)}{y^{\frac 32}}, \,\,g_2(t,y)= e^{-\frac{( \xi(t)-\xi(t-y))^2}{4y}}$$
defined in $C_T$.
Then

\be\label{g_lip1}|g(t_2)-g(t_1)|=\left |\int_{t_1}^{t_2}g_1(t_2,y)g_2(t_2,y)f(t_2-y)\,dy+\right.\ee
$$\left.\int_0^{t_1}g_1(t_2,y)g_2(t_2,y)f(t_2-y)-g_1(t_1,y)g_2(t_1,y)f(t_1-y)\right|.$$
Observe that, by hypothesis, there exist a constant $H_1$ such that 
$$|g_1(t_2,y)g_2(t_2,y)f(t_2-y)-g_1(t_1,y)g_2(t_1,y)f(t_1-y)|\le \frac{H_1}{\sqrt y}$$
then 
\be\label{g_lip2}\left|\int_{t_1}^{t_2}g_1(t_2,y)g_2(t_2,y)f(t_2-y)\,dy\right|\le H_2\sqrt {t_2-t_1},\ee
where $H_2$ is a proper constant.

\noindent Let us consider in $(0,t_1)$ the term 
$$g_1(t_2,y)g_2(t_2,y)f(t_2-y)-g_1(t_1,y)g_2(t_1,y)f(t_1-y)=$$
$$g_1(t_1,y)g_2(t_1,y)(f(t_2-y)-f(t_1-y))+
(g_1(t_2,y)g_2(t_2,y)-g_1(t_1,y)g_2(t_1,y))f(t_2-y).$$
We can choose $H_3$ such that
$$
|g_1(t_1,y)g_2(t_1,y)(f(t_2-y)-f(t_1-y))|\le \frac{H_3}{\sqrt y}\sqrt{t_2-t_1},
$$
then we can choose $H_4$, that does not depend on $t_1$, such that
\be\label{g_lip3}
\left| \int_0^{t_1}g_1(t_1,y)g_2(t_1,y)(f(t_2-y)-f(t_1-y))\,dy\right|\le \int_0^{t_1}\frac{H_1}{\sqrt y}\sqrt{t_2-t_1}\,dy\le H_4\sqrt{t_2-t_1}.
\ee
It remains to estimate  the term
$$
(g_1(t_2,y)g_2(t_2,y)-g_1(t_1,y)g_2(t_1,y))f(t_2-y)=$$
$$
[(g_1(t_2,y)-g_1(t_1,y)]g_2(t_2,y))f(t_2-y)+[(g_2(t_2,y)-g_2(t_1,y)]g_1(t_1,y))f(t_2-y).
$$
First of all we consider 
$g_1(t_2,y)-g_1(t_1,y)$. Let us apply, for every $y\in (0,t_1)$, the Lagrange theorem to the function $s(\cdot,t)=\frac{\xi(\cdot)-\xi(\cdot-y)}{y^{\frac 32}}$. Then, for every $y\in (0,t_1)$, there exists $t_y\in (t_1,t_2)$  such that
$$
|g_1(t_2,y)-g_1(t_1,y)|=\left|\frac{\xi'(t^y)-\xi'(t^y-y)}{y^{\frac 32}}\right |t_2-t_1|\le \frac{H_1}{y^{\frac32-\beta}}|t_2-t_1|
$$
where again  $H_5$ is  a proper constant.
Since $\beta>\frac12$ we deduce that 
there exists $H_6$ such that
\be\label{g_lip4}
\left |\int_0^{t_1}[g_1(t_2,y)-g_1(t_1,y)]g_2(t_2,y))f(t_2-y)\,dy\right|\le H_6|t_2-t_1|.
\ee
Finally we have to consider the term
$g_2(t_2,y)-g_2(t_1,y)$. Again we apply the Lagrange theorem to the function $l(\cdot,y)=e^\frac{-(\xi(\cdot)-\xi(\cdot-y))^2}{4y}$, obtaining
$$
|g_2(t_2,y)-g_2(t_1,y)|=$$
$$
\left|\displaystyle{ e^\frac{-(\xi(t'_y)-\xi(t'_y-y))^2}{4y}}\frac{(\xi(t'_y)-\xi(t'_y-y))(\xi'(t'_y)-\xi'(t'_y-y))}{2y}\right||t_2-t_1|,
$$
for every $y\in (0,t_1)$ we have $t'_y\in (t_1,t_2)$. 
Then we deduce 
\be\label{g_lip5}
\int_0^{t_1}[(g_2(t_2,y)-g_2(t_1,y)]g_1(t_1,y))f(t_2-y)\,dy\le H_7|t_2-t_1|
\ee
with  a proper constant $H_7$.

Obviously the constant $H_1,\dots, H_7$ depends only on the Holder constants of $\xi'$, $f$ and the Lipschitz constant of $\xi$.
Therefore, using \eqref{g_lip}--\eqref{g_lip5}, we obtain the thesis.
\qed
\begin{lemma}\label{lemma_kilbas2} Let $\xi\in C^{1+\beta}([0,T])$, $\beta>\frac{1}{2}$, $f\in C^{\frac12}([0,T])$. Let $K$ the function defined in \eqref{KeG} and $G(t,s)$ the associated function obtained in Proposition \ref{prop_volterra} and Remark \ref{G_alfa}.
Then the function 
$$
v(t):=\int_0^t G(t,s)f(s)\,ds,
$$
is in the space  $C^{\frac 12}([0,T]).$
\end{lemma}
\proof We easily have this result since, as observed $\xi\in C^{1+\beta}$, $\beta>\frac12$ implies that $G(\cdot,s)$ is Holder continuous with exponent $\frac 12$ uniformly respect to the variable $s$. 
\qed

 {\it Proof of Theorem \ref{prop_kilbas}.}\, In order to prove this result we have to separate the terms that depends on $m'$ from which that do not depend on it.
Therefore it is useful to rewrite $N_1=N^1_1+N^2_1$ and $N_2=N^1_2+N^2_2$ where $N_i^1$, $i=1,2$ depends only on the initial data $u_0'$ and $u_T'$. Since ${g_0^-}'=\frac{m'(t)}{\gamma_0}$ and ${g_0^+}'=\frac{m'(t)}{\gamma_2}$, we get 
$$N_1^1(t)=\frac{-\sqrt{\pi}}{\sqrt{t}}\int_0^{\infty}e^{-\frac{(\underline \xi(t)-x)^2}{4t}}u_T'(K-x)\,dx\, \hbox{ in }(0,|\gamma_0| T),$$
$$N_2^1(t)=\frac{\sqrt{\pi}}{\sqrt{t}}\int_0^{\infty}e^{-\frac{(\overline \xi(t)-x)^2}{4t}}u_0'(x)\,dx\, \hbox{ in }(0,\gamma_2 T)$$
and

$$N_1^2(t)=-\pi\int_0^t \frac{e^{-\frac{(\underline \xi(t)-\underline \xi(s))^2}{4(t-s)}}}{(t-s)^{\frac12}}\frac{{m}'(T-\frac{s}{|\gamma_0|})}{|\gamma_0|^2}\,ds\, \hbox{ in }(0,|\gamma_0| T),$$
$$N_2^2(t)=-\pi\int_0^t \frac{e^{-\frac{(\overline \xi(t)-\overline \xi(s))^2}{4(t-s)}}}{(t-s)^{\frac12}}\frac{{m}'(\frac{s}{\gamma_2})}
{\gamma_2^2}\,ds\, \hbox{ in }(0,\gamma_2 T)$$
where $\overline \xi$ and $\underline\xi$ are defined in Proposition \ref{DN2_map} and \ref{DN1_map}.

\noindent Let us introduce  the following twelve  functions 
$$s_1(t)=\frac {|\gamma_0|}{\pi}N_1^1(|\gamma_0|(T- t)),\,\,s_2(t)=-\frac {\gamma_2}{\pi}N_2^1(\gamma_2 t),$$
  $$
s_3(t)=\frac{|\gamma_0|}{\pi^2}\int_0^{|\gamma_0|(T- t)} K_1(|\gamma_0|(T- t),s)N_1^1(s)\,ds, $$
$$
s_4(t)=-\frac{\gamma_2}{\pi^2}\int_0^{\gamma_2 t} K_2(\gamma_2 t,s)N_2^1(s)\,ds, $$

$$ s_5(t)=\frac{|\gamma_0|}{\pi^2}\int_0^{|\gamma_0|(T- t)} G_1(|\gamma_0|(T- t),s)N_1^1(s)\,ds$$

$$ s_6(t)=-\frac{\gamma_2}{\pi^2}\int_0^{\gamma_2 t} G_2(\gamma_2 t,s)N_2^1(s)\,ds,$$
$$s_7(t)=\frac {|\gamma_0|}{\pi}N_1^2(|\gamma_0|(T- t)),\,\,s_8(t)=-\frac {\gamma_2}{\pi}N_2^2(\gamma_2 t),$$
 $$
s_9(t)=\frac{|\gamma_0|}{\pi^2}\int_0^{|\gamma_0|(T- t)} K_1(|\gamma_0|(T- t),s)N_1^2(s)\,ds, $$
$$
s_{10}(t)=-\frac{\gamma_2}{\pi^2}\int_0^{\gamma_2 t} K_2(\gamma_2 t,s)N_2^2(s)\,ds, $$

$$ s_{11}(t)=\frac{|\gamma_0|}{\pi^2}\int_0^{|\gamma_0|(T- t)} G_1(|\gamma_0|(T- t),s)N_1^2(s)\,ds,$$
$$ s_{12}(t)=-\frac{\gamma_2}{\pi^2}\int_0^{\gamma_2 t} G_2(\gamma_2 t,s)N_2^2(s)\,ds.$$

Therefore the right hand side of \eqref{afraid} is equal to $\displaystyle{\sum_{i=1}^{12}s_i(t)}$  in $[0,T]$. Moreover the functions $s_i$, $i=1,\cdots,6$ do not depend on the function $m(t)$.

\noindent Let us analyze the terms $s_1$ and $s_2$.  It is clear that these have the same regularity of $N_1^1(t)$ and $N_2^1(t)$. By changing variable  we can write 
$$N_2^1(t)=\displaystyle{C_2\int_{\frac{-\xi(\frac t{\gamma_2})}{2\sqrt t}}^{+\infty} e^{-y^2}u'_0\left(2\sqrt t y+\xi\left(\frac t{\gamma_2}\right)\right)\,dy}$$
where $C_2$ is a proper constant. Then, we easily obtain, using $u''_0\in L^{\infty}((0,\infty))$ and $\xi(0)=0$, the following estimate  
$$\frac{d}{dt}N_2^1(t)\le\frac{\overline C_2}{\sqrt t} \hbox{ in } (0,\gamma_2 T]$$
where $\overline C_2$ is a proper constant. This 
implies  that $N_2^1(t)$ is  in the space $C^{\frac 12}([0,\gamma_2 T])$.

 Analogously  with the change of variable $y=\frac{x+\xi(T-\frac{t}{|\gamma_0|})-K}{2\sqrt t}$ we obtain
$$
N_2^1(t)=C_1\int_{\frac{\xi(T-\frac{t}{|\gamma_0|})-K}{2\sqrt t}}^{+\infty}e^{-y^2}u'_T\left(-2\sqrt t y+\xi\left(T-\frac t{|\gamma_0|}\right)\right)\,dy$$
and using the conditions $\xi(T)=K$, $u''_T\in L^{\infty}((-\infty,K))$ we obtain that $N_1^1$ is in the space $C^{\frac 12}([0,|\gamma_0|T])$. 

\noindent In order to prove that the functions $s_3$ and $s_4$ are in $C^{\frac12}([0,T])$ it is enough to apply Lemma \ref{lemma_kilbas}. The Holder continuity of the functions $s_5$ and $s_6$ are consequence of Lemma \ref{lemma_kilbas2}.

Let us analyze the terms $s_7,\dots s_{12}$ involving the unknown function $m(t)$.

\noindent We have
$$s_7(t)=-|\gamma_0|\int_0^{|\gamma_0|(T-t)}\frac{e^{\frac{-(\underline \xi(|\gamma_0|(T-t))-\underline \xi(s))^2}{4[|\gamma_0|(T-t)-s]}}}{(|\gamma_0|(T-t)-s)^{\frac12}}\frac{m'(T-\frac{s}{|\gamma_0|})}{|\gamma_0|^2}\,ds.$$
By passing to the  variable $y=T-\frac{s}{|\gamma_0|}$ and using the relation $\underline \xi(|\gamma_0|(T-z))=K-\xi(z)$
we obtain the following rewriting
\be \label{ktt1}
s_7(t)=-\int_t^T \frac{e^{\frac{-( \xi(y)- \xi(t))^2}{4[|\gamma_0|(y-t)]}}}{|\gamma_0|^{\frac 12}(y-t)^{\frac12}}m'(y)\,dy.
\ee
The term $s_8$ is given by a similar singular integral kernel. In fact we get
$$
s_8(t)=\int_0^{\gamma_2t}\frac{e^{\frac{-(\overline \xi(\gamma_2t)-\overline \xi(s))^2}{4(\gamma_2t-s]}}}{(\gamma_2t-s)^{\frac12}}\frac{m'(\frac{s}{\gamma_2})}{\gamma_2}\,ds.
$$  
Again by changing variable and using the definition of $\overline \xi$ this becomes
\be\label{ktt2}s_8(t)=\int_0^t \frac{e^{\frac{-( \xi(y)- \xi(t))^2}{4[\gamma_2y]}}}{\gamma_2^{\frac 12}(t-y)^{\frac12}}m'(y)\,dy.
\ee
Let us consider term $s_{10}(t)$.
We have
$$
s_{10}(t)=\frac1{\pi}\int_0^{\gamma_2t}K_2(\gamma_2t,s)\left[\int_0^s\frac{e^{\frac{-(\overline \xi(s)-\overline \xi(\tau))^2}{4(s-\tau)}}}{(s-\tau)^{\frac12}}\frac{m'(\frac{\tau}{\gamma_2})}{\gamma_2}\,d\tau\right]\,ds
$$
Let us use first the change of variable $s=\gamma y$, then $\tau=\gamma_2z$. We obtain
$$
s_{10}(t)=\frac{\gamma_2^{\frac 12}}{\pi}\int_0^{t}K_2(\gamma_2t,\gamma_2 y)\left[\int_0^s\frac{e^{\frac{-(\overline \xi(\gamma_2 y)-\overline \xi(\gamma_2 z))^2}{4\gamma_2(y-z)}}}{(y-z)^{\frac12}}m'(z)\,dz\right]\,dy
$$
Utilizing the definition of $\overline \xi$ and changing the order of integration,  $s_{10}$ becomes
\be\label{compli_1}\frac{{\gamma_2}^{\frac12}}{\pi}\int_0^{t}\left[\int_z^t K_2(\gamma_2t,\gamma_2 y)\frac{e^{\frac{-( \xi(y)-\xi(z))^2}{4\gamma_2(y-z)}}}{(y-z)^{\frac12}}\,dy\right]m'(z)\,dz
\ee
denoting with $\tilde K_2(t,z)$ the function in the square brackets of \eqref{compli_1} we observe that this is an Holder  continuous function defined in  $C_T$ by Lemma \ref{lemma_kilbas0}. 

\noindent With analogues technics we easily prove that there exist three   Holder continuous functions of order $\frac 12$; $\tilde K_1$, $\tilde G_1$, defined in $[0,T]\times [0,T]\setminus  C_T$ and  
$\tilde G_2$ defined in $\overline C_T$ such that 
\be\label{ktt3}
s_9(t)=\int_t^T \tilde K_1(t,z)m'(z)\,dz,\, s_{11}(t)=\int_t^T \tilde G_1(t,z)m'(z)\,dz,
\ee
\be\label{ktt4}
s_{12}(t)=\int^t_0 \tilde G_2(t,z)m'(z)\,dz.
\ee
Considering equation \eqref{afraid},  we can write $m(t)=m(0)+\int_0^t m'(z)\,dz$. Therefore
the assertion of the Theorem comes from  the Holder regularity of the functions $\xi'$, $s_1,\dots,s_6$ and from the formulas \eqref {ktt1}--\eqref{ktt4}.
\qed

\begin{remark}\label{diag} Observe that from the computations given in the proof of Theorem \ref{prop_kilbas} we obtain $k_1(t,t)=-\frac 1{|\gamma_0|^{\frac 12}}$ and $k_2(t,t)=\frac 1{\gamma_2^{\frac 12}}$. More precisely these values are obtained in \eqref{ktt1}-\eqref{ktt2}.
\end{remark}
 
Generalized  Abel' s equation like  \eqref{abel1} are strictly related to the solutions of singular integral equations. We do not analyze in details this kind of problems we refer to the book of Samko, Kilbas, Marichev,   \cite{kilbas} for a complete  treatise of the subject and to \cite{Mu} and \cite{Gh} for the singular integral equations. We just recall that choosing properly the class of the solutions  it can be proved the Noether nature of the operator  
\be\label{abel2}
A(\psi)=\int_0^t\frac{k_1(t,s)}{\sqrt{t-s}}{\psi(s)}\,ds+\int_t^T\frac{k_2(t,s)}{\sqrt{s-t}}{\psi(s)}\,ds.
\ee
 In particular in \cite{kilbas} it is proved that if $k_i$ are Holder continuos of order $\gamma>\frac 12$, denoting  with 
 $$I^{\frac12}(L^p(0,T))=\left\{v\in L^p(0,T)\,:\, \hbox{ there exists }\phi\in L^p(0,T) \hbox{ such that } \right.$$
 $$\left.v(t)=\int_0^t\frac{\phi(s)}{\sqrt{t-s}}\,ds\right\}
 $$
 then $A:L^p(0,T)\rightarrow I^{\frac12}(L^p(0,T))$
 is a Noether operator when $p<2$ (see also \cite{Ru} for a generalization in the case $p\ge 2$). 
 
 \noindent Unfortunately in our situation $k_i$ are Holder continuos at most of order $\frac 12$ then we can not apply this result.  On the other hand in this case we can apply a result proved again in \cite {kilbas} which gives the ``algebraic Noether nature'' of equation \eqref{abel1}. In order to enunciate this result we use Remark \ref{diag} and rewrite equation \eqref{abel1} as
 \be\label{abella}
 \frac{1}{\sqrt{\gamma_2}}\int_0^t\frac{m'(s)}{\sqrt{t-s}}\,ds-\frac{1}{\sqrt{|\gamma_0|}}\int_t^T\frac{m'(s)}{\sqrt{s-t}}\,ds+\int_0^T T(t,s)m'(s)\,ds=h(t)
 \ee
 where $T$ is a regular function in $\{(t,s)\in [0,T]\times[0,T]\,:\, t\neq s\}$ which has a singularity along the diagonal of order strictly lesser that $\frac{1}{ |t-s|^\frac12}$.
 In particular by straightforward calculation it satisfies the hypothesis of Lemma 31.3 in \cite{kilbas}.   
 
 \noindent In this case equation \eqref{abella} is well posed when $m'$ is taken in the space $H^*$ and $h$ is in the space $H^*_\frac12$,
 where 
$$H^*_\frac12:=\large\bigcup_{\frac12<\lambda\le 1,\epsilon_1,\epsilon_2\in (0,1)} H_0^{\lambda}(\epsilon_1,\epsilon_2),
$$
$$
H_0^{\lambda}(\epsilon_1,\epsilon_2):=\left\{f(t)=t^{1-\epsilon_1}(T-t)^{1-\epsilon_2}g(t)\,:\, g\in C^{\lambda}([0,T]), g(0)=g(T)=0 \right\}$$
and
$$
H^*:=\large\bigcup_{0<\lambda\le 1,\epsilon_1,\epsilon_2\in (0,1)} H_0^{\lambda}(\epsilon_1,\epsilon_2).
$$
In some sense the space $H^*_\frac12$ correspond to the space of Holder continuous function of order strictly bigger that $\frac 12$ a part on the end points $t=0$ and $t=T$  where it can also  singular. Analogously the space $H^*$ correspond to the space of Holder continuous function of any arbitrary order  a part on the end points $t=0$ and $t=T$, where singularity is allowed. These are good spaces for studying Abel's equations since it can be proved (see \cite {kilbas}) that operators
$$A_1(\psi)=\int_0^t\frac{\psi(s)}{\sqrt{t-s}}\,ds,\,\, A_2(\psi)=\int_t^T\frac{\psi(s)}{\sqrt{s-t}}\,ds,
$$
are well defined from $H^*$ to $H^*_\frac12$.

\noindent In our case we easily obtain that $h(t)$ is in the space $H^*_\frac12$. In fact, for hypothesis $\xi'\in C^{\beta}([0,T])$, with $\beta>\frac 12$. Moreover we obtained in the proof of Theorem \ref{prop_kilbas} that $N_1^1$, $N_1^2$ are not only Holder continuous function of order $\frac 12$ but also derivable with derivate that has a singularity of order  $\frac 12$ at the end point $t=0$ and $t=T$. Then we obtain that $s_1,\cdots s_6$ are in the space  $H^*_\frac12 $.

\noindent Therefore we can apply the following result obtained in \cite {kilbas} (p. 650, Theorem 31.11)
\begin{theorem}\label{bella} Equation \eqref{abella} is solvable in $H^*$ if and only if 
$$
\int _0^T h(s)\psi_j(s)\,ds=0,\, j=1,\cdots, k
$$
where $\psi_j$ is a complete system of solution of the homogenous equation
\be\label{abella_hom}
 \frac{1}{\sqrt{\gamma_2}}\int_t^T\frac{m'(s)}{\sqrt{t-s}}\,ds-\frac{1}{\sqrt{|\gamma_0|}}\int_0^t\frac{m'(s)}{\sqrt{s-t}}\,ds+\int_0^T T(t,s)m'(s)\,ds=0,
 \ee
 where $\psi$ has the following form
 $$\psi(t)=\frac{\psi^*(t)}{t^{\frac12}(T-t)^{\frac12}}$$
 with $\psi^*(t)$ an Holder continuous function of any order in $[0,T]$, $k$ is the finite dimension of the subspace.

\noindent  The difference between numbers of linearly independent solutions of equations \eqref{abella} and \eqref{abella_hom} is equal to 1.
\end{theorem}
\begin{remark} Observe that, choosing properly $u_0$ and $u_T$,  previous result   allow to obtain explicit entropy solution of the forward--backward parabolic equation in which  the presence of the unstable phase is non--trivial. It is important to point out that the solution that we obtain is also an entropy solution of the original problem \eqref{B_T} only if $v(\xi(t),t)\in (A,B)$ for every $t\in [0,T]$. When this solution exists Theorem \ref{bella} gives not uniqueness. In general we can choose $\tau\le T$, such that  $v(\xi(t),t)\in (A,B)$ for every $t\in [0,\tau]$, in order to have an entropy solution of the forward--backward equation in the strip $\R\times [0,\tau]$. 
\end{remark}

\begin{remark} The case $\xi(\cdot)\equiv 0$ is much more easier since we can choose in \eqref{abella} $T(t,s)\equiv 0$. In this context  we can obtain more satisfactory  results for the generalized Abel's equation (see \cite{kilbas}).
\end{remark}

\begin{remark} The analysis handled in this section could be useful to study forward parabolic equation with discontinuous coefficients along a given interface.
\end{remark}

\address{Dipartimento di Matematica ``G. Castelnuovo''\\
Universit\`a di Roma ``La Sapienza''\\
Piazzale A. Moro, 2 --- 00185 Roma, Italy\\

\email{terracin@mat.uniroma1.it}\\

\end{document}